\documentclass[12pt]{article}
\usepackage[dvips]{epsfig}
\usepackage{graphicx}
\usepackage{amsmath}
\newcommand{\stack}[2]{\genfrac{}{}{0pt}{}{#1}{#2}}
\usepackage{amssymb}
\parskip=\medskipamount

\begin{document}

\def\pmb#1{\setbox0=\hbox{#1}%
  \kern-.015em\copy0\kern-\wd0
  \kern.05em\copy0\kern-\wd0
  \kern-.015em\raise.0433em\box0 }

\def\bdJ{\pmb{$J$}}  
\def\bdw{\pmb{$w$}}  
\def\bds{\pmb{$s$}}  
\def\bdF{\pmb{$F$}}  
\def\bd0{\pmb{$0$}}  
\def\bdf{\pmb{$f$}}
\def\bdy{\pmb{$y$}}
\def\bdx{\pmb{$x$}}
\def\bdz{\pmb{$z$}}
\def\bdu{\pmb{$u$}}
\def\bdA{\pmb{$A$}}
\def\bda{\pmb{$a$}}
\def\caL{\cal L}
\def\R{\mathbb R}
\def\C{\mathbb C}
\def\bdab{\bar{\pmb{$a$}}}
\def\bdom{\pmb{$\omega$}}
\def\be{\begin{equation}}
\def\ee{\end{equation}}
\def\bea{\begin{eqnarray}}
\def\eea{\end{eqnarray}}
\def\beas{\begin{eqnarray*}}
\def\arccos{{arccos}}
\def\eeas{\end{eqnarray*}}
\def\ds{\displaystyle}
\def\b{\phantom{-}}
\newtheorem{dfn}{Definition}
\newtheorem{thm}{Theorem}
\newtheorem{lem}{Lemma}
\newtheorem{rem}{Remark}
\newtheorem{cor}[thm]{Corollary}
\date{}
\addtolength{\baselineskip}{0.3\baselineskip}
\renewcommand{\theequation}{\arabic{section}.\arabic{equation}}
\title{\Large\bf Two-step hybrid methods adapted to the numerical integration of perturbed oscillators}
\author{Hans Van de Vyver\\
{Department of Mathematics, Katholieke Universiteit Leuven,}\\
{Celestijnenlaan 200 B, B-3001 Heverlee, Belgium}\\
e-mail: hans$\_$vandevyver@hotmail.com}

\maketitle

\begin{abstract}
Two-step hybrid methods specially adapted to the numerical integration of perturbed oscillators are obtained. The formulation of the methods is based on a refinement of classical Taylor expansions due to Scheifele [{\em Z. Angew. Math. Phys.}, {\bf 22}, 186--210 (1971)]. The key property is that those algorithms are able to integrate exactly harmonic oscillators with frequency $\omega$ and that, for perturbed oscillators, the local error contains the (small) perturbation parameter as a factor. The methods depend on a parameter $\nu=\omega\,h$, where $h$ is the stepsize. Based on the B2-series theory of Coleman [{\em IMA J. Numer. Anal.}, {\bf 23}, 197--220 (2003)] we derive the order conditions of this new type of methods. The linear stability and phase properties are examined. The theory is illustrated with some fourth- and fifth-order explicit schemes. Numerical results carried out on an assortment of test problems (such as the integration of the orbital motion of earth satellites) show the relevance of the theory.
\end{abstract}

\vspace*{1.5cm}
\noindent AMS Classification~: 65L05 

\vspace*{0.5cm}

\vspace*{0.5cm}

\newpage

\setcounter{equation}{0}
\section{Introduction}
In the last decades, there has been a great interest in the research of methods for the numerical integration of initial value problems (IVP) associated to second-order ordinary differential equations (ODE)
\be
y''=f(x,y),\qquad y(x_0)=y_0,\qquad y'(x_0)=y'_0,
\label{ode2}
\ee
in which the first derivative does not appear explicitly. These problems appear often in practice. Of course, since (\ref{ode2}) can be written as an IVP for a system of two equations of first-order, the problem can be solved by algorithms for first-order equations. However, this will be less efficient than if methods specially devised for the given problem would be used. The construction of methods specialized for (\ref{ode2}) is a well established area of investigation. Many multistep methods (such as St\"ormer--Cowell methods) and two-step hybrid methods for (\ref{ode2}) have been developed, see for example Lambert \& Watson (1976), Chawla (1984), Chawla \& Rao~(1987), Coleman~(1989), Simos~(1999), Tsitouras~(2003), Coleman~(2003) and Franco~(2006a) to mention a few. Two-step hybrid methods are considered to be more efficient than the rival Runge-Kutta-Nystr\"om methods for (\ref{ode2}). For example, the standard fourth-order explicit Runge-Kutta-Nystr\"om method (see Hairer {\em et al.} (1993)) requires three function evaluations whereas the fourth-order explicit Numerov method of Chawla (1984) requires only two function evaluations per step.

Quite often the solution of (\ref{ode2}) exhibits an oscillatory behaviour; think, for instance, of the pendulum problem in celestial mechanics or of the Schr\"odinger equation in quantum mechanics. For problems having highly oscillatory solutions standard methods with unspecialized use can require a huge number of steps to track the oscillations. One way to obtain a more efficient integration process is to construct numerical methods with an increased algebraic order. On the other hand, the construction and implementation of high algebraic order methods is not evident. Alternatively, one can consider methods that use the detailed information of the high-frequency oscillation. There is a vast literature on this subject; an extensive bibliography is summarized by Petzold {\em et al.} (1997). Scheifele (1971) was concerned with the solution of {\em perturbed oscillators}, i.e., second-order problems of the form
\be
y''=-\,\omega^2\,y+g(x,y),\qquad y(x_0)=y_0,\qquad y'(x_0)=y'_0,
\label{pert}
\ee
where the magnitude of the perturbation force satisfies $|g(x,y)|<<\omega^2\,|y|$.  Scheifele rewrote the solution of (\ref{pert}) as a series of a set of functions, the $G$-functions, more adequate to perturbed oscillators than the classical polynomial Taylor expansion. The Scheifele $G$-functions method is capable to integrate exactly the harmonic oscillator or unperturbed problem (i.e. (\ref{pert}) with $g=0$). In spite of its excellent behaviour, the Scheifele $G$-functions method has the disadvantage that it is strictly application-dependent. Several authors have applied Scheifele's approach for constructing numerical methods adapted to perturbed oscillators. Most of these papers are focused on space dynamical problems such as an accurate integration of orbit problems or long-term prediction of satellite orbits. Some Scheifele $G$-functions based multistep codes are designed by Martin \& Ferr\'andiz~(1997). Also adapted methods without first derivatives have been constructed by L\'opez {\em et al.}~(1999). A first Runge-Kutta type version of the Scheifele $G$-functions method is due to Gonz\'alez {\em et al.} (1999). A theoretical foundation for these adapted Runge-Kutta-Nystr\"om (ARKN) methods is given by Franco~(2002,2005,2006b).

Our objective in this paper is apply Scheifele's approach to two-step hybrid methods. This was already proposed by Van de Vyver (2007a) for the simple explicit Numerov method. The excellent numerical results reported in that paper strongly suggest to construct higher-order methods of this type. This is possible when a more theoretical framework would be developed. This is the purpose of this work. The paper is organized as follows. Section~\ref{sectsh} is of an introductory nature: we recall the class of classical two-step hybrid (TSH) methods. In Section~\ref{secatsh} we recall Scheifele's approach. This idea will be extended to TSH methods, the resulting methods are denoted by ATSH methods. Section~\ref{secorder} is devoted to the order conditions for ATSH methods. This part heavily relies on the work of Coleman~(2003) for classical TSH methods. Some general stability results for ATSH methods are reported in Section~\ref{secstab}. The concepts of such a stability analysis find its origin in the work of Coleman \& Ixaru (1996) and Franco (2005). Section~\ref{secpldis} provides general results on the phase properties of ATSH methods. The analysis is based on the work of Franco~(2005). Section~\ref{seccon} deals with the construction of fourth- and fifth-order explicit ATSH methods. Several possibilities are explored such as minimizing the error constant, increasing the phase-lag order, dissipative or not,~\ldots\,The classical companions of the new methods are previously derived by Franco~(2006a). Section~\ref{secnumexp} collects numerical examples for a variety of problems chosen to illustrate particular features of the ATSH methods obtained. The new methods are compared with other high-quality methods. The paper concludes with a brief summary of the work considered here.
\section{Classical two-step hybrid methods}
\label{sectsh}
Two-step hybrid (TSH) methods for (\ref{ode2}) are defined by
\bea
Y_i&=&(1+c_i)\,y_n-c_i\,y_{n-1}+h^2\,\ds\sum_{j=1}^s a_{ij}\,f(x_n+c_j\,h,Y_j),\qquad i=1,\ldots,s,\label{2stepint}\\
y_{n+1}&=&2\,y_n-y_{n-1}+h^2\,\ds\sum_{i=1}^s b_i\,f(x_n+c_i\,h,Y_i),\label{2stepfin}
\eea
where $y_{n-1}$, $y_n$ and $y_{n+1}$ are approximations of $y(x_n-h)$, $y(x_n)$ and $y(x_n+h)$, respectively. TSH methods can be in short-hand notation represented by the Butcher table
\[
\begin{array}{c|ccc}
c_1& a_{11}& \ldots & a_{1s}\\
\vdots & \vdots & \ddots &\vdots\\
c_s& a_{s1}& \ldots & a_{ss}\\
\hline
& b_1 & \ldots & b_s
\end{array}
\qquad
=
\qquad
\begin{array}{c|c}
c & A\\
\hline\\[-4mm]
& b^T
\end{array},
\]
where $c,b\in\R^{s\times 1}$ and $A\in\R^{s\times s}$. These coefficients are derived by imposing the necessary and sufficient conditions for {\em convergence}, i.e. {\em consistency} and {\em zero-stability}, see Henrici~(1962) for the general theory.

For exact starting values, the {\em local truncation error~(lte)} of the method at $x_n$ is
\be
lte=y(x_n+h)-2\,y(x_n)+y(x_n-h)-h^2\,\ds\sum_{i=1}^sb_i\,f(x_n+c_i\,h,Y_i).
\label{lte}
\ee
The method is of {\em algebraic order} $p$ if $lte={\mathcal O}(h^{p+2})$. The {\em principal local truncation error (plte)} is the leading term of (\ref{lte}). For a $p$th-order method this is of the form
\be
plte=\ds\frac{h^{p+2}}{(p+2)!}\,\ds\underset{\stack{t\in T_2}{\rho(t)=p+2}}{\sum}\alpha(t)\,\left(1+(-1)^{p+2}-b^T\,\Psi(t)\right)F(t)(y_n,y_n'),
\label{plte}
\ee
where $\alpha(t)$, $\rho(t)$, $\Psi''(t)$, $F(t)$ and $T_2$ are defined in Coleman~(2003). The coefficients of $F(t)(y_n,y_n')$ in (\ref{plte}) will be denoted as $e_{p+1}(t)$. The quantity 
\be
E_{p+1}=\left(\ds\underset{\stack{t\in T_2}{\rho(t)=p+2}}{\sum}e^2_{p+1}(t)\right)^{1/2},
\label{errc}
\ee
will be called the {\em error constant} of the $p$th-order method. Traditionally, the order conditions for TSH methods are usually derived by expansions in Taylor series. These expansions are calculated essentially by brute force. On the other hand, Coleman~(2003) obtained the order conditions for TSH methods by using the theory of B-series. Analogously to the case of RK(N) methods, the determination of the order of a TSH method is based on checking certain relationships between the coefficients of the method.

The linear stability analysis of methods for solving~(\ref{ode2}) is based on the scalar test equation (see Lambert \& Watson (1976))
\be
y''=-\lambda^2\,y,\qquad \lambda>0.
\label{testeq}
\ee
An application of a TSH method to~(\ref{testeq}) yields
\be
\begin{array}{lll}
Y&=&(e+c)\,y_n-c\,y_{n-1}-H^2\,A\,Y,\qquad H=\lambda\,h,\\[3mm]
y_{n+1}&=&2\,y_n-y_{n-1}-H^2\,b^T\,Y,
\end{array}
\label{2steptesteq}
\ee
where $Y=(Y_1,\ldots,Y_s)^T$ and $e=(1,\ldots,1)^T\in\R^{s\times1}$.
Elimination of the vector $Y$ from~(\ref{2steptesteq}) results in the difference equation
\be
y_{n+1}-S(H^2)\,y_n+P(H^2)\,y_{n-1}=0,
\label{diffeq}
\ee
where
\be
\begin{array}{lll}
S(H^2)&=&2-H^2\,b^T\,(I+H^2\,A)^{-1}\,(e+c),\\[3mm]
P(H^2)&=&1-H^2\,b^T\,(I+H^2\,A)^{-1}\,c.\\
\end{array}
\label{SP}
\ee
The solution of the difference equation (\ref{diffeq}) is determined by the {\em characteristic equation}
\be
\xi^2-S(H^2)\,\xi+P(H^2)=0.
\label{Char}
\ee
Of particular interest for periodic motion is the situation where the roots of (\ref{Char}) lie on the unit circle. For example, in celestial mechanics it is desired that numerical orbits do not spiral inwards or outwards. This periodicity condition is equivalent to 
\be
P(H^2)=1 \qquad \mbox{and} \qquad |S(H^2)|<2,\qquad \forall H\in(0,H_{per}^2),
\label{percond}
\ee
and the interval $(0,H_{per}^2)$ is called the {\em interval of periodicity}.  If the necessary condition $P(H^2)=1$ to have of a non-empty interval of periodicity is not satisfied, we can ask when the numerical solution remains bounded. This stability condition is equivalent to
\[
P(H^2)<1 \qquad \mbox{and} \qquad |S(H^2)|<1+P(H^2),\qquad \forall H\in(0,H_{stab}^2),
\]
and the interval $(0,H_{stab}^2)$ is called the {\em interval of absolute stability}. 

Another related concept, which is important when solving problems of the form~(\ref{ode2}) is the phase-lag of the method. In phase analysis one compares the phases of $\exp(\pm\,i\,H)$ with the phases of the roots of the characteristic equation~(\ref{Char}). Following the approach of van der Houwen \& Sommeijer~(1987) for RKN methods, the quantities
\be
\phi(H)=H-\arccos\left(\ds\frac{S(H^2)}{2\,\sqrt{P(H^2)}}\right), \qquad d(H)=1-\sqrt{P(H^2)},
\label{phase}
\ee
are the {\em phase-lag (or dispersion)} and the {\em dissipation (or amplification error)}, respectively. The method is said to have {\em phase-lag order} $q$ and {\em dissipation order} $r$~if
\[
\phi(H)=c_{\phi}\,H^{q+1}+{\mathcal O}(H^{q+3}),\qquad d(H)=c_d\,H^{r+1}+{\mathcal O}(H^{r+3}).
\]
The constants $c_{\phi}$ and $c_d$ are called the {\em phase-lag} and {\em dissipation constants}, respectively.
Methods with $d(H)=0$ are {\em zero-dissipative}.
\section{Two-step hybrid methods for perturbed oscillators}
\label{secatsh}
\subsection{Notations and exact solution}
Although, Scheifele's method is based on $G$-functions, in this paper we consider the related $\phi$-functions which are suggested by Franco (2002) for the derivation of the order conditions for ARKN methods. The coefficients of Scheifele's $G$-functions method are dependent on the frequency $\omega$ and stepsize $h$. By using the $\phi$-functions, the coefficients are dependent on only one variable $\nu=\omega\,h$. 

The solution of (\ref{pert}) can be expressed as
\be
y(x_n+h)=y(x_n)\,\cos(\nu)+hy'(x_n)\,\ds\frac{\sin(\nu)}{\nu}+\ds\frac{1}{\omega}\,\ds\int_{x_n}^{x_{n+1}}g(x,y(x))\,\sin(\omega\,(x_{n+1}-x))\,dx.
\label{inteq}
\ee
We carry out the change of variable $x=x_n+h\,z$ in (\ref{inteq}) and we denote $\varphi(x)=g(x,y(x))$. Now the exact solution becomes
\be
y(x_n+h)=y(x_n)\,\cos(\nu)+hy'(x_n)\,\ds\frac{\sin(\nu)}{\nu}+h^2\,\ds\int_0^1\varphi(x_n+h\,z)\,\ds\frac{\sin(\nu\,(1-z))}{\nu}\,dz.
\label{eqvarphi}
\ee
Suppose that the function $\varphi(x)$ admits an expansion of the form
\be
\varphi(x_n+h\,z)=\ds\sum_{j=0}^{\infty}h^j\,\varphi^{(j)}(x_n)\,\ds\frac{z^j}{j!}.
\label{series}
\ee
We can write that
\be
y(x_n+h)=y(x_n)\,\cos(\nu)+hy'(x_n)\,\ds\frac{\sin(\nu)}{\nu}+\ds\sum_{j=0}^{\infty}h^{j+2}\varphi^{(j)}(x_n)\,\int_0^1\ds\frac{\sin(\nu\,(1-z))}{\nu}\,\ds\frac{z^j}{j!}\,dz.
\ee
Introducing the following notations 
\be
\phi_0(\nu)=\cos(\nu),\qquad \phi_1(\nu)=\ds\frac{\sin(\nu)}{\nu},\qquad \phi_{j+2}(\nu)=\int_0^1\ds\frac{\sin(\nu\,(1-z))}{\nu}\,\ds\frac{z^j}{j!}\,dz,\qquad j\geq 0,
\ee
we arrive to the expression of the exact solution of the perturbed problem (\ref{pert}) in terms of $\phi$-functions
\be
y(x_n+h)=y_n\,\phi_0(\nu)+h\,y'_n\,\phi_1(\nu)+\ds\sum_{j=0}^{\infty}h^{j+2}\,\varphi^{(j)}(x_n)\,\phi_{j+2}(\nu).
\label{exsol}
\ee
Remark that the analytical solution of the harmonic oscillator is approximated exactly by the expansion (\ref{exsol}).

Some interesting properties of the $\phi$-functions are listed in the following theorem.
\begin{thm}
\begin{enumerate}
\item $\ds\lim_{\nu\rightarrow 0}\phi_j(\nu)=\ds\frac{1}{j!},\qquad j\geq 0$.
\item The $\phi$-functions can be expressed as
\be
\phi_{2j}(\nu)=\ds\frac{(-1)^j}{\nu^{2\,j}}\Biggl(\cos(\nu)-\ds\sum_{k=0}^{j-1}(-1)^k\,\ds\frac{\nu^{2\,k}}{(2\,k)!}\Biggr),\qquad j\geq 0,
\ee
\be
\phi_{2j+1}(\nu)=\ds\frac{(-1)^j}{\nu^{2\,j+1}}\Biggl(\sin(\nu)-\ds\sum_{k=0}^{j-1}(-1)^k\,\ds\frac{\nu^{2\,k+1}}{(2\,k+1)!}\Biggr),\qquad j\geq 0.
\ee
\item The Taylor series expansions of the $\phi$-functions are
\be
\phi_j(\nu)=\ds\sum_{k=0}^{\infty}(-1)^k\ds\frac{\nu^{2\,k}}{(2\,k+j)!},\qquad j\geq 0.
\ee
\item $\phi_{j+1}(\nu)=\ds\int_0^1\cos(\nu\,(1-z))\,\ds\frac{z^j}{j!}\,dz,\qquad j\geq 0.$
\item We have the following recurrence relation
\be
\phi_j(\nu)+\nu^2\,\phi_{j+2}(\nu)=\ds\frac{1}{j!},\qquad j\geq 0.
\label{recurrence}
\ee
\end{enumerate}
\label{thmphi}
\end{thm}
The $\phi$-functions are related to the Scheifele $G$-functions by $G_j(h)=h^j\,\phi_j(\nu)$, $j\geq 0$. For further details and proofs about $G$-functions, see Scheifele (1971), Fair\'en~{\em et al.}~(1994) and Mart\'in \& Ferr\'andiz (1997).

According to Theorem \ref{thmphi} (point 1) it is clear that when the frequency $\omega\rightarrow 0$ ($\nu\rightarrow 0$) the series (\ref{exsol}) will become
\be
y(x_n+h)=y(x_n)+h\,y'(x_n)+\ds\sum^{\infty}_{j=0}\ds\frac{h^{j+2}}{(j+2)!}y^{(j+2)}(x_n),
\ee
which is the classical Taylor expansion of the exact solution. Thus Scheifele's series~(\ref{exsol}) is a refinement of the classical Taylor method.

\subsection{Formulation of the method}
An $s$-stage TSH method (\ref{2stepint})--(\ref{2stepfin}) can be rewritten in the following alternative form  
\[
\begin{array}{rcl}
k_i'&=&f\Bigl(x_n+c_i\,h,(1+c_i)\,y_n-c_i\,y_{n-1}+h^2\,\ds\sum_{j=1}^sa_{ij}\,k'_j\Bigr),\qquad i=1,\ldots,s,\\[3mm]
y_{n+1}&=&2\,y_n-y_{n-1}+h^2\,\ds\sum_{i=1}^s k_i'.
\end{array}
\]
We can see that $k'_i$ are evaluations of the function $f$ at the points $x_n+c_i\,h$, where the second argument is an approximation to the solution at this point. Then, we have
\[
y(x_n+c_i\,h)\approx (1+c_i)\,y_n-c_i\,y_{n-1}+h^2\,\ds\sum_{j=1}^sa_{ij}\,k'_j,\qquad i=1,\ldots,s.
\]
For perturbed oscillators, i.e. when $f(x,y)=-\omega^2\,y+g(x,y)$, the internal stages can be approximated by
\be
k_i=g(x_n+c_i\,h,Y_i),\qquad i=1,\ldots,s,
\label{keq}
\ee
where
\[
Y_i=(1+c_i)\,y_n-c_i\,y_{n-1}+h^2\,\ds\sum_{j=1}^sa_{ij}\,(-\omega^2\,Y_j+k_j).
\]
The coefficients $a_{ij}$ represent the weights of the quadrature formulas used in the approximation of the internal stages.

The final stage is determined as follows. We can avoid the calculation of the first derivative of the solution of (\ref{eqvarphi}) by adding this expression with positive and negative stepsize to get
\be
y(x_n+h)=2\,\phi_0(\nu)\,y(x_n)-y(x_n-h)+h^2\,\int^1_{-1}\ds\frac{\sin(\nu\,(1-|z|))}{\nu}\,\varphi(x_n+h\,z)\,d\,z.
\label{even} 
\ee
We shall approximate the exact solution by using the quadrature formula
\[
\int^1_{-1}\ds\frac{\sin(\nu\,(1-|z|))}{\nu}\,\varphi(x_n+h\,z)\,d\,z\approx \ds\sum_{i=1}^s b_i\,k_i,
\]
where the $k$-values are given by (\ref{keq}). 

Altogether, we arrive to the following definition.
\begin{dfn}
An $s$-stage adapted two-step hybrid (ATSH) method for the numerical integration of the IVP (\ref{pert}) is given by the scheme
\be
\begin{array}{l}
Y_i=(1+c_i)\,y_n-c_i\,y_{n-1}+h^2\,\ds\sum_{j=1}^sa_{ij}\,\Bigl(-\omega^2\,Y_j+g(x_n+c_j\,h,Y_j)\Bigr),\qquad 1,\ldots,s,\\[3mm]
y_{n+1}=2\,\phi_0(\nu)\,y_n-y_{n-1}+h^2\,\ds\sum_{i=1}^s b_i\,g(x_n+c_i\,h,Y_i),
\end{array}
\label{atsh}
\ee
which can be expressed in Butcher notation by the table of coefficients
\[
\begin{array}{c|ccc}
c_1& a_{11}& \ldots & a_{1s}\\
\vdots & \vdots & \ddots &\vdots\\
c_s& a_{s1}& \ldots & a_{ss}\\
\hline
& b_1 & \ldots & b_s
\end{array}
\qquad
=
\qquad
\begin{array}{c|c}
c & A\\
\hline\\[-4mm]
& b^T
\end{array}.
\]
\end{dfn}
Remark that when $\omega\rightarrow 0$, ATSH methods reduce to classical TSH methods. 

As said, the convergence of a method is covered by consistency and zero-stability. The consistency (i.e. algebraic order is at least 1) follows form Section~\ref{secorder}. The theorem of Ixaru \& Rizea (1987) says that any method applied to $y''=0$ with resulting difference equation 
\[
y_{n+1}+a_1(h)\,y_n+y_{n-1}=0, 
\]
is zero-stable if $a_1(h)=-2+{\mathcal O}(h^{q})$, $q>2$. Using Theorem~\ref{thmphi} (point~3) it is easy to see that ATSH methods are zero-stable.
\section{Order conditions for ATSH methods}
\label{secorder}
Similarly to the classical case, the {\em principal local truncation error (lte)} of the ATSH method (\ref{atsh}) is given by
\[
lte=y(x_n+h)-2\,\phi_0(\nu)\,y(x_n)+y(x_n-h)-h^2\,\ds\sum_{i=1}b_i\,g(x_n+c_i\,h,Y_i).
\]
The method is of {\em algebraic order} $p$ if $lte={\mathcal O}(h^{p+2})$. Our next aim is to derive order conditions for ATSH methods by adapting the recently developed B2-series theory of Coleman~(2003). In what follows, the reader is referred to that paper for all the definitions and notations. The theory of B2-series is applicable only to one-step methods so we have to search for a one-step formulation of ATSH methods. A modification of Coleman's proofs at several places will deliver the requested order conditions.
\subsection{Adapted B2-series}
Repeated differentiation of  $\varphi$ with respect to the independent variable $x$ gives
\[
\begin{array}{l}
\varphi^{(0)}=g(y),\\[3mm]
\varphi^{(1)}=g^{(1)}(y)(y'),\\[3mm]
\varphi^{(2)}=g^{(2)}(y)(y',y')+g^{(1)}(y)(f(y)),\\[3mm]
\varphi^{(3)}=g^{(3)}(y)(y',y',y')+3\,g^{(2)}(y)(y',f(y))+g^{(1)}(y)(f^{(1)}(y)(y')),\\[3mm]
\hspace{7mm}\ldots
\end{array}
\]
The difference with the classical theory lies in the fact that every elementary differential starts with a Fr\'echet-derivative of $g$ instead of $f$. The following definition explains how each elementary differential can be associated with a rooted tree.
\begin{dfn}The function $G$ on $T_2\backslash\{\O,\tau'\}$ is defined by
\begin{enumerate}
\item $G(\tau)(y,y')=g$.
\item If $t=[t_1,\ldots,t_m]_2\in T_2$, then
\[
G(t)(y,y')=g^{(m)}(y)\left(F(t_1)(y,y'),\ldots,F(t_m)(y,y')\right),
\]
where the function $F$ is recursively defined in Definition~3 of Coleman~(2003).
\end{enumerate}
\end{dfn}
Analogously to the classical theory, it is obvious that
\be
\varphi^{(j)}=\ds\underset{\stack{t\in T_2}{\rho(t)=j+2}}{\sum}\alpha(t)\,G(t)(y,y'),
\label{deriv}
\ee
where $\alpha(t)$ represents the number of distinct monotonic labellings of the vertices of $t\in T_2$.

B2-series are defined in Definition~4 of Coleman (2003). Here that definition is adopted more pertinent for our methods.
\begin{dfn}
Let $\beta$ a mapping from $T_2$ to $\R$. The adapted B2-series with coefficient function $\beta$ is a formal series of the form
\[
\tilde{B}(\beta,y)=\ds\sum_{t\in T_2\backslash\{\O,\tau'\}}\ds\frac{h^{\rho(t)}}{\rho(t)!}\alpha(t)\,\beta(t)\,G(t)(y,y').
\]
\end{dfn}
Coleman's fundamental lemma is then reformulated for the adapted case as follows.
\begin{lem}
Let $B(\beta,y)$ be a classical B2-series. Then $h^2\,g(B(\beta,y))$ is an adapted B2-series,
\[
h^2\,g(B(\beta,y))=\tilde{B}(\beta'',y),
\]
with
\[
\beta''(\O)=\beta''(\tau')=0,\qquad \beta''(\tau)=2,
\]
and for all other $t=[t_1,\ldots,t_m]_2\in T_2$,
\[
\beta''(t)=\rho(t)\,(\rho(t)-1)\,\ds\prod_{i=1}^m\beta(t_i).
\]
\label{lemfund}
\end{lem}
The proof is essentially the same as the original proof.
\subsection{One-step formulation}
By defining $F_n:=(y_{n+1}-\phi_0(\nu)\,y_{n})/h$ the second equation of (\ref{atsh}) can be expressed as a pair of equations
\[
\begin{array}{lll}
y_n&=&\phi_0(\nu)\,y_{n-1}+h\,F_{n-1},\\[3mm]
F_n&=&\phi_0(\nu)\,F_{n-1}-\omega\,\nu\phi_1^2(\nu)y_{n-1}+h\,(b^T\otimes I)\,g(Y).
\end{array}
\]
Now, the one-step formulation takes the form
\be
u_n=M(\nu)\,u_{n-1}+h\Phi(u_{n-1},h),
\label{1stepa}
\ee
with
\be
\begin{array}{l}
M(\nu)=\left(\begin{array}{cc}\phi_0(\nu)&0\\-\omega\,\nu\,\phi^2_1(\nu)&\phi_0(\nu)\end{array}\right),\\[7mm]
u_n=\left(\begin{array}{c}y_n\\F_n\end{array}\right)\qquad\mbox{and}\qquad \Phi(u_{n-1},h)=\left(\begin{array}{c}F_{n-1}\\( b^T\otimes I)\,g(Y)\end{array}\right),
\end{array}
\label{1stepb}
\ee
and $Y$ is defined implicitly by
\be
\begin{array}{lll}
Y&=&(e+c)\otimes y_n-c\otimes y_{n-1}+h^2\,(A\otimes I)\,(-\omega^2\,Y+g(Y))\\[3mm]
&=&(e+(\phi_0(\nu)-1)\,c)\otimes y_{n-1}+h\,(e+c)\otimes F_{n-1}+h^2\,(A\otimes I)\,(-\omega^2\,Y+g(Y)).\\
\end{array}
\label{1stepc}
\ee
\subsection{Order conditions}
The vector $u_n$ is an approximation for $z_n=z(x_n,h)$, where
\be
z(x,h)=\left(
\begin{array}{c}
y(x)\\[3mm]
\ds\frac{y(x+h)-\phi_0(\nu)\,y(x)}{h}
\end{array}
\right).
\label{eqz}
\ee
For exact starting values, the $lte$ of the one-step formulation (\ref{1stepa})--(\ref{1stepc}) is
\be
d_n=z_n-M(\nu)\,z_{n-1}-h\,\Phi(z_{n-1},h),
\label{eqd}
\ee
with
\be
\Phi(z_{n-1},h)=\left(\begin{array}{c}\ds\frac{y(x_n)-\phi_0(\nu)\,y(x_{n-1})}{h}\\[3mm](b^T\otimes I)\,g(Y)\end{array}\right),
\label{eqphi}
\ee
where Y is now defined implicitly by
\[
Y=e\otimes\,y(x_{n-1})+(e+c)\otimes \bigl(y(x_n)-y(x_{n-1})\bigr)+h^2\,(A\otimes I)\,(-\omega^2\,Y+g(Y)).
\]
\begin{dfn}
The ATSH method~(\ref{atsh}) is of algebraic order $p$ when $d_n={\mathcal O}(h^{p+1})$.
\end{dfn}
We are now ready to present one of the main results of this paper.
\begin{thm}The ATSH method~(\ref{atsh}) is of algebraic order $p$ if and only if, for trees $t\in T_2$
\[
b^T\,\psi''(t)=\Bigl(1+(-1)^{\rho(t)}\Bigr)\,\rho(t)!\,\phi_{\rho(t)}(\nu),
\]
for $\rho(t)\leq p+1$ but not for some trees of order $p+2$.
\label{thmorder}
\end{thm}
{\em Proof. }Observing (\ref{eqz})--(\ref{eqphi}) we have that the first component of $d_n$ is zero. 
Each component of the vector $Y$ can be expanded as a B2-series
\be
Y_i(x_n)=B\Bigr(\psi_i,y(x_n)\Bigl)=\ds\sum_{t\in T_2}\ds\frac{h^{\rho(t)}}{\rho(t)!}\,\alpha(t)\,\psi_i(t)\,F(t)(y_n,y'_n).
\label{B2Y}
\ee
The coefficients $\psi_i(t)$ can be generated recursively by formulas (3.6)--(3.7) of~Coleman~(2003). We substitute the B2-series~(\ref{B2Y}) into the second component of $d_n$ and we apply Lemma~\ref{lemfund}. An easy calculation gives
\[
\ds\frac{1}{h}\,\left(y(x_n+h)-2\,\phi_0(\nu)\,y(x_n)+y(x_n-h)-h^2\,\ds\sum_{i=1}^s b_i\,g(Y_i(x_n))\right)
\]
\be
=\ds\frac{1}{h}\,\left(2\,\ds\sum_{j=1}^{\infty} h^{2\,j}\varphi_n^{(2\,j-2)}\phi_{2\,j}(\nu)-\ds\sum_{i=1}^sb_i\,\tilde B\Bigr(\psi_i'',y(x_n)\Bigl)\right).
\label{onestepeq}
\ee
With (\ref{deriv}) in mind, the left side of (\ref{onestepeq}) becomes
\be
2\,\ds\sum_{j=1}^{\infty} h^{2\,j}\varphi_n^{(2\,j-2)}\phi_{2\,j}(\nu)=\ds\underset{\stack{t\in T_2}{\rho(t)\,\mbox\small{even}}}{\sum}\,h^{\rho(t)}\,\alpha(t)\,\phi_{\rho(t)}(\nu)\,G(t)(y_n,y_n').
\label{order1}
\ee
The right side of (\ref{onestepeq}) may be written as
\be
\ds\sum_{i=1}^s b_i\,\tilde B(\psi_i'',y_n)=\ds\sum_{t\in T_2}\ds\frac{h^{\rho(t)}}{\rho(t)!}\,\alpha(t)\,b_i\,\psi_i''(t)\,G(t)(y_n,y_n').
\label{order2}
\ee
The theorem follows when comparing (\ref{order1}) and (\ref{order2}).

The order conditions up to order six are listed in Table~\ref{ordercol}.
\begin{rem}
Reconsidering Section~5 of Coleman~(2003) it is obvious that, in order to reduce the number of order conditions, the simplifying conditions for ATSH methods are the same as for classical TSH methods.

\end{rem}
\begin{table}[ht!]
\begin{center}
\begin{tabular}{|c|c|l|}
\hline
\mbox{Tree }$t$&$\rho(t)$&\mbox{Order condition}\\
\hline
&&\\[-3mm]
$t_{21}$&2&$\sum_i b_i=2\,\phi_2(\nu)$\\[3mm]
\hline
&&\\[-3mm]
$t_{31}$&3&$\sum_i b_i\,c_i=0$\\[3mm]
\hline
&&\\[-3mm]
$t_{41}$&4&$\sum_i b_i\,c_i^2=4\,\phi_4(\nu)$\\[3mm]
$t_{42}$& & $\sum_{i,j} b_i\,a_{ij}=2\,\phi_4(\nu)$\\[3mm]
\hline
&&\\[-3mm]
$t_{51}$&5& $\sum_i b_i\,c_i^3=0$\\[3mm]
$t_{52}$& & $\sum_{i,j} b_i\,c_i\,a_{ij}=2\,\phi_4(\nu)$\\[3mm]
$t_{53}$& & $\sum_{i,j} b_i\,a_{ij}\,c_j=0$\\[3mm]
\hline
&&\\[-3mm]
$t_{61}$&6& $\sum_i b_i\,c_i^4=48\,\phi_6(\nu)$\\[3mm]
$t_{62}$& & $\sum_{i,j} b_i\,c_i^2\,a_{ij}=24\,\phi_6(\nu)$\\[3mm]
$t_{63}$& & $\sum_{i,j} b_i\,c_i\,a_{ij}\,c_j=-\frac{2}{3}\,\phi_4(\nu)+8\,\phi_6(\nu)$\\[3mm]
$t_{64}$& & $\sum_{i,j,k} b_i\,a_{ij}\,a_{ik}=\phi_4(\nu)+12\,\phi_6(\nu)$\\[3mm]
$t_{65}$& & $\sum_{i,j} b_i\,a_{ij}\,c_j^2=4\,\phi_6(\nu)$\\[3mm]
$t_{66}$& & $\sum_{i,j,k} b_i\,a_{ij}\,a_{jk}=2\,\phi_6(\nu)$\\[3mm]
\hline
&&\\[-3mm]
$t_{71}$&7& $\sum_i b_i\,c_i^5=0$\\[3mm]
$t_{72}$& & $\sum_{i,j} b_i\,c_i^3\,a_{ij}=24\,\phi_6(\nu)$\\[3mm]
$t_{73}$& & $\sum_{i,j} b_i\,c_i^2\,a_{ij}\,c_j=0$\\[3mm]
$t_{74}$& & $\sum_{i,j,k} b_i\,c_i\,a_{ij}\,a_{ik}=24\,\phi_6(\nu)$\\[3mm]
$t_{75}$& & $\sum_{i,j,k} b_i\,c_i\,a_{ij}\,a_{jk}=-\frac{1}{6}\,\phi_4(\nu)+4\,\phi_6(\nu)$\\[3mm]
$t_{76}$& & $\sum_{i,j} b_i\,c_i\,a_{ij}\,c_j^2=\frac{1}{3}\,\phi_4(\nu)$\\[3mm]
$t_{77}$& & $\sum_{i,j,k} b_i\,a_{ij}\,a_{ik}\,c_k=-\frac{1}{3}\,\phi_4(\nu)+4\,\phi_6(\nu)$\\[3mm]
$t_{78}$& & $\sum_{i,j} b_i\,a_{ij}\,c_j^3=0$\\[3mm]
$t_{79}$& & $\sum_{i,j,k} b_i\,a_{ij}\,c_j\,a_{jk}=2\,\phi_6(\nu)$\\[3mm]
$t_{7,10}$& & $\sum_{i,j,k} b_i\,a_{ij}\,a_{jk}\,c_k=0$\\[3mm]
\hline
\end{tabular}
\caption{\label{ordercol}Order conditions}
\end{center}
\end{table}

\subsection{Error analysis}
From the proof of Theorem \ref{thmorder} it follows that the $plte$ of a $p$th-order ATSH method is given by
\[
plte^{ATSH}=\ds\frac{h^{p+2}}{(p+2)!}\,\ds\underset{\stack{t\in T_2}{\rho(t)=p+2}}{\sum}\alpha(t)\,\left(1+(-1)^{p+2}-{b^{(0)}}^T\,\Psi^{(0)}(t)\right)G(t)(y_n,y_n'),
\]
where ${b^{(0)}}^T$ and $\Psi^{(0)}$ represents the $b^T$- and $\Psi$-values of the corresponding classical TSH method.
The $plte$ of this classical method for (\ref{ode2}) reads
\be
plte^{TSH}=\ds\frac{h^{p+2}}{(p+2)!}\,\ds\underset{\stack{t\in T_2}{\rho(t)=p+2}}{\sum}\alpha(t)\,\left(1+(-1)^{p+2}-{b^{(0)}}^T\,\Psi^{(0)}(t)\right)F(t)(y_n,y_n').
\label{plte2}
\ee
In order to obtain a connection between $plte^{TSH}$ and $plte^{ATSH}$ we need a relationship between $F(t)$ and $G(t)$. This can be easily seen as follows. We consider trees in which the root starts with a chain of 3 vertices (including the root) having exactly one son. We call such a tree a {\em semi-tall tree}. We denote by $T_2^*$ the set of semi-tall trees. The {\em truncated tree} $t^-$ of a semi-tall tree $t$ is obtained by deleting the first two vertices. Clearly, the number of semi-tall trees of order $p+2$ is equal to the number of trees of order $p$. Using the above terminology, it is easy to see that
\[
G(t)(y,y')=
\left\{
\begin{array}{ll}
F(t)(y,y')+\omega^2\,F(t^-)(y,y') & \mbox{if}\,\,\, t\in T_2^*,\\[5mm]
F(t)(y,y')                        & \mbox{if}\,\,\, t\notin T_2^*.\\
\end{array}
\right.
\]
We conclude with
\be
plte^{ATSH}=plte^{TSH}+\omega^2\,\ds\frac{h^{p+2}}{(p+2)!}\,\ds\underset{\stack{t\in T_2^*}{\rho(t)=p+2}}{\sum}\alpha(t)\,\left(1+(-1)^{p+2}-{b^{(0)}}^T\,{\Psi}^{(0)}(t)\right)F(t^-)(y_n,y_n').
\label{pltead}
\ee
For the calculation of the error constant, $E^{ATSH}_{p+1}$, we have to consider the coefficients of $F(t)(y_n,y_n')$ and the coefficients of $\omega^2\,F(t^-)(y_n,y_n')$ in (\ref{pltead}). Observing (\ref{plte2})--(\ref{pltead}) it is clear that
\be
E^{ATSH}_{p+1}=\left(\ds\underset{\stack{t\in T_2}{\rho(t)=p+2}}{\sum}{l_i\,(e^{TSH}_{p+1}})^2(t_i)\right)^{1/2}\qquad\mbox{with}\qquad 
l_i=
\left\{
\begin{array}{ll}
2 & \mbox{if}\,\,\, t\in T_2^*,\\[3mm]
1 & \mbox{if}\,\,\, t\notin T_2^*.\\
\end{array}
\right.
\label{errca}
\ee
\section{Linear stability analysis}
\label{secstab}
Linear stability and phase-lag analysis of ATSH methods is also based on the model equation~(\ref{testeq}). However, this equation has to be rewritten in the following appropriate form 
\be
y''=-\omega^2\,y-\epsilon\,y,\qquad \omega^2+\epsilon>0,
\label{testperteq}
\ee
where $\omega$ represents an estimation of the dominant frequency $\lambda$ of (\ref{testeq}), and $\epsilon=\lambda^2-\omega^2$ is the error of that estimation. This modified test equation is prompted by the work of Franco (2005) for ARKN methods. At the first sight, one should believe that the estimated frequency $\omega$ should be equal to dominant frequency $\lambda$. This is generally a satisfying approach but in practical applications it is possible to obtain more accurate results for different values of $\lambda$ and $\omega$. 
The cubic oscillator
\[
y''=-y+\epsilon\,y^3,\qquad y(0)=1,\qquad y'(0)=1,
\]
provides such an example. Although this is a nonlinear problem, for small $\epsilon$-values we may apply linear stability analysis, resulting in $\lambda=1$. However, Vigo-Aguiar {\em et al.} (2004) have proved that more accurate results are obtained when selecting $\omega=\sqrt{1-0.75\,\epsilon}$.

An ATSH method (\ref{atsh}) applied to (\ref{testperteq}) yields
\[
\begin{array}{lll}
Y&=&(e+c)\,y_n-c\,y_n-(\nu^2+z)\,A\,Y,\\[5mm]
y_{n+1}&=&2\,\phi_0(\nu)\,y_n-y_{n-1}-z\,b^T\,Y,\qquad \nu=\omega\,h,\qquad z=\epsilon\,h^2.
\end{array} 
\]
Elimination of the vector $Y$ gives the recurrence relation
\be
y_{n+1}-S(\nu^2,z)\,y_n+P(\nu^2,z)\,y_{n-1}=0,
\label{rec}
\ee
where
\be
S(\nu^2,z)=2\,\phi_0(\nu)-z\,b^T\,N^{-1}\,(e+c),\qquad P(\nu^2,z)=1-z\,b^T\,N^{-1}\,c,
\label{eqSP}
\ee
and
\be
N=I+(\nu^2+z)\,A,\qquad e=(1,\ldots,1)^T. 
\label{eqN}
\ee
The characteristic equation is
\be
\xi^2-S(\nu^2,z)\,\xi+P(\nu^2,z)=0.
\label{Chareq}
\ee
Firstly, let us consider dissipative ATSH methods. Working with (\ref{Chareq}), we can ask, for a given method (i.e., a given $\omega$), and a given test frequency $\lambda$, what restriction must be placed on the stepsize $h$ to ensure that the stability condition
\be
P(\nu^2,z) < 1\qquad\mbox{and}\qquad |S(\nu^2,z)| < P(\nu^2,z)+1,
\label{stabcond}
\ee 
is satisfied. This question can be answered by examining $S(\nu^2,z)$ and $P(\nu^2,z)$ in the $\nu-z$ plane. For ARKN methods such a stability analysis was introduced by Franco~(2005). The following definition is originally formulated by Coleman~\&~Ixaru~(1996) for exponentially fitted methods for (\ref{ode2}). Here, it is adjusted in terms of the methods of concern.
\begin{dfn}
For a dissipative ATSH method with $S(\nu^2,z)$ and $P(\nu^2,z)$ where $\nu=\omega\,h$ and $z=\epsilon\,h$, and $\omega$ and $\epsilon$ are given, the primary interval of absolute stability is the largest interval $(0,h_0)$ such that (\ref{stabcond}) holds for all stepsizes $h\in(0,h_0)$. If, when $h_0$ is finite, (\ref{stabcond}) holds also for $\gamma<h<\delta$, where $\gamma>h_0$ then the interval $(\gamma,\delta)$ is a secondary interval of absolute stability. The region of absolute stability is a region in the $\nu-z$ plane ($\nu>0$), throughout which (\ref{stabcond}) holds. Any closed curve defined by
\[
P(\nu^2,z) = 1\qquad{or}\qquad |S(\nu^2,z)| = P(\nu^2,z)+1,
\]
is a stability boundary.
\end{dfn}
Likewise, for zero-dissipative ATSH methods the definition of the {\em primary interval of periodicity} and the {\em region of periodicity} is evident.

In the particular case when the main frequency is exactly known (i.e. $z=0$) we have for both dissipative and zero-dissipative methods that
\[
S(\nu^2,0)=2\,\cos(\nu)\qquad\mbox{and} \qquad P(\nu^2,0)=1.
\]
It follows that the $\nu$-axis is a stability boundary. On this line the periodicity condition~(\ref{percond}) is satisfied except when $\nu=n\,\pi$ for positive integer $n$.

In the dissipative case, when the frequency is not exactly known the stepsize has to be selected carefully. Here we show some sensible points. 
\begin{thm}For dissipative ATSH methods there exist values for $\omega$ and $\epsilon$ for which the primary interval of absolute stability is empty.
\end{thm}
{\em Proof. }Consider the function $F$ defined as
\[
F(H^2)=b^T\,(I+H^2\,A)^{-1}\,c.
\]
Assume that $F$ is continuous at $H^2=\nu^2$. So we can find an interval $(-z_0,z_0)$ such that $F(\nu^2+z)$ has the same sign for all $z\in(-z_0,z_0)$. It turns out that for such $z$-values the function $P$, as given in (\ref{eqSP})--(\ref{eqN}), has a different sign at the points $(\nu,-z)$ and $(\nu,z)$. From the absolute stability condition (\ref{stabcond}) it follows  that an ATSH method which is stable at $(\nu,-z)$, is not stable at $(\nu,z)$. Thus the $\nu$-axis acts as a stability boundary in the sense that it separates stable and unstable regions. This concludes the proof.
\section{Phase-lag and dissipation analysis}
\label{secpldis}
For any method corresponding to the characteristic equation~(\ref{Chareq}), the quantities
\be
\phi(\nu^2,z)=H-\arccos\left(\ds\frac{S(\nu^2,z)}{2\,\sqrt{P(\nu^2,z)}}\right), \qquad d(\nu^2,z)=1-\sqrt{P(\nu^2,z)},
\label{phasead}
\ee
are called the phase-lag and the amplification error, respectively. As pointed out by Franco (2005) for ARKN methods, the analysis of the phase-lag and the dissipation becomes more useful if we introduce
\be
\nu=\ds\frac{\omega}{\sqrt{\omega^2+\epsilon}}\,H,\qquad z=\ds\frac{\epsilon}{\omega^2+\epsilon}\,H^2,
\label{subs}
\ee
in (\ref{phasead}). So we arrive to the following definition.
\begin{dfn}
The phase-lag order is $q$ if 
\be
\phi(\nu^2,z)=c_{\phi}(\omega^2,\epsilon)\,H^{q+1}+{\mathcal O}(H^{q+3}),
\label{phlf}
\ee
and the dissipation order is $r$ if
\be
d(\nu^2,z)=c_d(\omega^2,\epsilon)\,\nu^{r+1}+{\mathcal O}(H^{r+3}).
\label{disf}
\ee
$c_{\phi}(\omega^2,\epsilon)$ and $c_{d}(\omega^2,\epsilon)$ are called the phase-lag and dissipation functions, respectively.
\end{dfn}
In the particular case when the main frequency is exactly known (i.e. $z=0$) the test equation (\ref{testperteq}) is integrated exactly and so there is no phase-error and no dissipation.

We investigate the phase properties when the main frequency is not exactly known. Let us define $C_j:=b^T\,A^{j-1}\,c$ and $U_j:=b^T\,A^{j-1}\,e$. Some algebraic manipulation gives
\begin{itemize}
\item ATSH method of algebraic order $p=2\,k$: 
\be
\begin{array}{lll}
S(\nu^2,z)&=&2\,\ds\sum_{j=0}^k\ds\frac{(-1)^j}{(2\,j)!}\,H^{2\,j}+2\,\ds\sum_{j=k+1}^{\infty}\ds\frac{(-1)^j}{(2\,j)!}\,H^{2\,j}\,\left(\ds\frac{\omega^2}{\omega^2+\epsilon}\right)^{j-k}\\[5mm]
&&+\ds\frac{\epsilon}{\omega^2+\epsilon}\,\ds\sum_{j=k+1}^{\infty}(-1)^j\,(U_j+C_j)\,H^{2\,j},\\[5mm]
P(\nu^2,z)&=&1+\ds\frac{\epsilon}{\omega^2+\epsilon}\,\ds\sum_{j=k+1}^{\infty}(-1)^j\,C_j\,H^{2\,j},
\end{array}
\label{speven}
\ee
\item ATSH method of algebraic order $p=2\,k-1$:
\be
\begin{array}{lll}
S(\nu^2,z)&=&2\,\ds\sum_{j=0}^k\ds\frac{(-1)^j}{(2\,j)!}\,H^{2\,j}+2\,\ds\sum_{j=k+1}^{\infty}\ds\frac{(-1)^j}{(2\,j)!}\,H^{2\,j}\,\left(\ds\frac{\omega^2}{\omega^2+\epsilon}\right)^{j-k}\\[5mm]
&&+\ds\frac{\epsilon}{\omega^2+\epsilon}\,(-1)^k\,C_k\,H^{p+1}+\ds\frac{\epsilon}{\omega^2+\epsilon}\,\ds\sum_{j=k+1}^{\infty}(-1)^j\,(U_j+C_j)\,H^{2\,j},\\[5mm]
P(\nu^2,z)&=&1+\ds\frac{\epsilon}{\omega^2+\epsilon}\,\ds\sum_{j=k}^{\infty}(-1)^j\,C_j\,H^{2\,j}.
\end{array}
\label{spodd}
\ee
\end{itemize}
When substituting (\ref{speven})--(\ref{spodd}) in (\ref{phasead}) and then considering the Taylor expansion with respect to $H$ it is sufficient to retain the term with the lowest power. After tedious but straightforward calculations we have concluded with
\begin{thm}
\begin{enumerate}
\item Assume that the algebraic order $p$ of a dissipative TSH method is even (odd) and that the phase-lag order is $q=p$ ($q=p+1$). Then the corresponding ATSH method has also phase-lag order~$q$. The leading term of the phase-lag (\ref{phlf}) is
\be
c_{\phi}(\omega,\epsilon)=\ds\frac{\epsilon}{\omega^2+\epsilon}\,c_{\phi},
\label{plf}
\ee
where $c_{\phi}$ is the phase-lag constant of the classical TSH method.
\item A dissipative TSH method and the corresponding ATSH method have both the same dissipation order. The leading term of the dissipation (\ref{disf}) is
\be
d_{\phi}(\omega,\epsilon)=\ds\frac{\epsilon}{\omega^2+\epsilon}\,d_{\phi},
\label{dis}
\ee
where $c_d$ is the dissipation constant of the classical TSH method.
\end{enumerate}
\end{thm}
From (\ref{plf}) it follows that the conditions for a ATSH method to have phase-lag order $q=p+2$ ($p$:~even) or $q=p+3$ ($p$: odd) are exactly the same as those of the corresponding classical method. This establishes
\begin{cor}
Assume that the algebraic order $p$ of a TSH method is even (odd) and that the phase-lag order is $q=p+2$ ($q=p+3$). Then the corresponding ATSH method has also phase-lag order~$q$.
\label{corr1}
\end{cor}
In general, Scheifele's adaptation does not conserve the phase-lag order for dissipative TSH methods. In contrast, we will show that the phase-lag order is always conserved in the zero-dissipative case. Taking into account the order conditions obtained in Section~\ref{secorder} and proceeding as in Section~9 of Coleman~(2003) we can reformulate Coleman's Theorem~6 for zero-dissipative ATSH methods as follows.
\begin{thm}For the determination of the the phase-lag order of a zero-dissipative ATSH method (\ref{atsh}) we have to compute the scalar quantities $C_k=b^T\,A^{k-1}\,c$ and $U_k=b^T\,A^{k-1}\,e$ for $k=1,2,\ldots$. The phase-lag order is q iff $U_k=2\,\phi_{2\,k}(\nu)$ for $k=1,\ldots,[\frac{p+1}{2}]$ and $C_k=0$ for $k=1,\ldots,[\frac{p}{2}]$ but one of those conditions is not satisfied when $p$ is replaced by $p+1$.
\label{phaseorderatsh}
\end{thm}

\begin{cor}
A zero-dissipative ATSH method and its classical companion have both the same phase-lag order.
\end{cor}
The phase-lag function is also of the form (\ref{plf}). 

Obviously we have in all cases that $c_{\phi}(\omega,0)=c_d(\omega,0)=0$, $c_{\phi}(0,\epsilon)=c_{\phi}$ and $c_d(0,\epsilon)=c_d$. When an acceptable estimate of the dominant frequency is available (i.e. $\epsilon\approx 0$) the magnitude of the phase-lag (\ref{plf}) and the amplification error (\ref{dis}) is then much smaller than those of the corresponding classical method. Furthermore, the more accurate the estimate of the dominant frequency, the smaller the phase-lag and the amplification error.

\section{Construction of explicit ATSH methods}
\label{seccon}
In this section we study the construction of explicit ATSH methods with algebraic orders four and five. Both dissipative and zero-dissipative methods are presented. The construction procedure in the classical case was previously considered by Franco~(2006a).

\subsection{Methods using two function evaluations per step}
Consider the explicit ATSH method defined by the table of coefficients
\[
\begin{array}{c|ccc}
-1&0&0&0\\
0&0&0&0\\
c_3&a_{31}&a_{32}&0\\
\hline
&b_1&b_2&b_3\\
\end{array}.
\]
Under the simplifying assumptions (see Coleman (2003))
\be
A\,e=\ds\frac{c^2+c}{2},
\label{simplcond}
\ee
the order conditions up to order four are
\be
b^T\,e=2\,\phi_2(\nu),\qquad b^T\,c=0,\qquad b^T\,c^2=4\,\phi_4(\nu),\qquad b^T\,c^3=0,\qquad b^T\,A\,c=0.
\label{ordercond4}
\ee
We have the unique solution
\be
b_1=b_3=2\,\phi_4(\nu),\qquad b_2=-4\,\phi_4(\nu)+2\,\phi_2(\nu),\qquad c_3=1,\qquad a_{31}=0,\qquad a_{32}=1.
\label{expnumad}
\ee
When $\nu\rightarrow 0$ the method reduces to the explicit Numerov method of Chawla~(1984). Remark that the values~(\ref{expnumad}) are obtained in a different way by Van de Vyver~(2007a). A stability and phase-lag analysis is also included in that paper.
\subsection{Methods using three function evaluations per step}
Next, we analyze the construction of explicit ATSH methods defined by the table of coefficients
\be
\begin{array}{c|cccc}
-1&0&0&0&0\\
0&0&0&0&0\\
c_3&a_{31}&a_{32}&0&0\\
c_4&a_{41}&a_{42}&a_{43}&0\\
\hline
&b_1&b_2&b_3&b_4\\
\end{array}.
\label{buts3}
\ee
\subsubsection{Dissipative fifth-order methods}
\label{secdis}
The order conditions up to order five are given by (\ref{simplcond})--(\ref{ordercond4}) with in addition
\be
b^T\,c^4=48\,\phi_6(\nu),\qquad b^T\,(c\,.\,A\,c)=-\frac{2}{3}\,\phi_4(\nu)+8\,\phi_6(\nu),\qquad b^T\,A\,c^2=4\,\phi_6(\nu).
\label{ordercond5}
\ee
Solving the equations (\ref{simplcond})--(\ref{ordercond4}) and (\ref{ordercond5}), the coefficients (\ref{buts3}) are determined in terms of the arbitrary parameter $c_3$. Two different strategies will be described  in order to get an optimal method. A first option is to determine $c_3$ so that the error constant~$E_6^{ATSH}$~(\ref{errca}) is as small as possible. The second option is to choose $c_3$ so that the method has phase-lag order eight.

\begin{figure}[htp!]
\begin{center}
\begin{tabular}{c}
\includegraphics[width=12 cm]{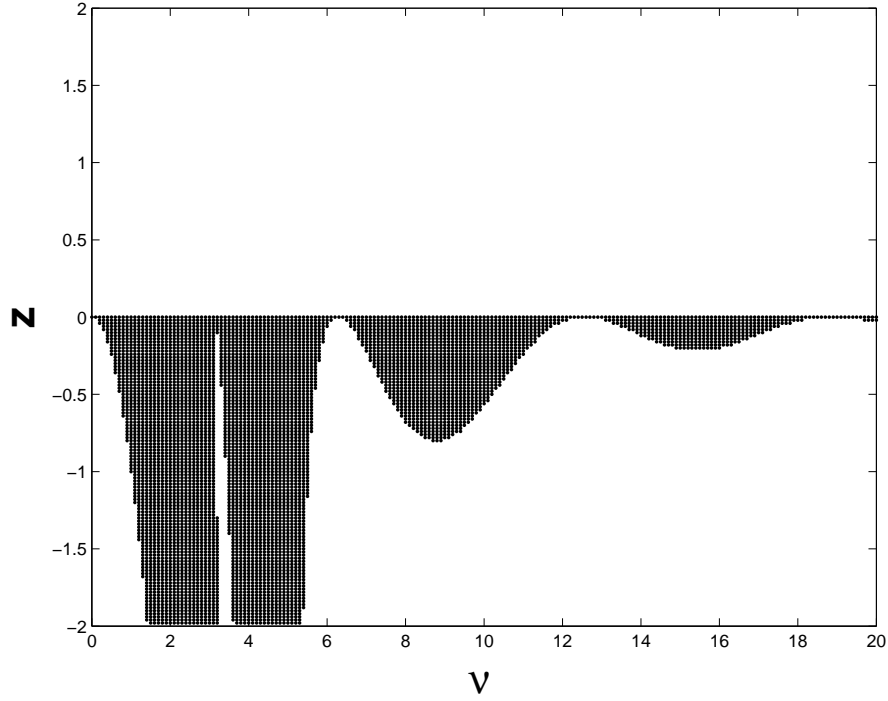}\\
\end{tabular}
\caption{\label{stab1}$\nu-z$ plot for ATSH method~(\ref{atsh1}).}
\end{center}
\end{figure}

{\bf * ATSH method with minimized error constant}
\newline
When minimizing the error constant $E_6^{ATSH}$, we obtain a value for $c_3$ which is very close (within a distance $<10^{-3}$) to those of a classical method of Franco~(2006a), $c_3=63/100$. For this reason we adopt Franco's method and we conclude with the coefficients
\be
\begin{array}{l}
a_{31}=\ds\frac{126651}{2000000},\hspace{1mm} a_{32}=\ds\frac{900249}{2000000},\hspace{1mm}  a_{41}=\ds\frac{100\,S_1\,S_2\,(720000\,\phi_6^2-124158\,\phi_6\,\phi_4+6031\,\phi_4^2)}{305488243\,\phi_4^4},\\[5mm]

a_{42}=\ds\frac{S_1\,S_2\,(-8000000\,\phi_6^2+886200\,\phi_6\,\phi_4+2849\,\phi_4^2)}{13119127\,\phi_4^4},\hspace{3mm} a_{43}=\ds\frac{20000\,S_1\,S_2\,S_3\,\phi_6}{2138417701\,\phi_4^4},\\[5mm]

b_1=\ds\frac{6\,(40000\,\phi_6-1323\,\phi_4)\,\phi_4}{163\,S_1},\\[5mm] b_2=\ds\frac{2\,(15338\,\phi_4^2-240000\,\phi_6\,\phi_4-3969\,\phi_4\,\phi_2+75600\,\phi_2\,\phi_6)}{189\,S_2},\\[5mm]

b_3=\ds\frac{400000000\,(12\,\phi_6-\phi_4)\,\phi_4}{30807\,S_3},\hspace{3mm} b_4=\ds\frac{3748322\,\phi_4^4}{9\,S_1\,S_2\,S_3},\qquad c_3=\ds\frac{6}{100}, \qquad c_4=\ds\frac{3\,S_2}{37\,\phi_4},\\[5mm]

S_1=600\,\phi_6-13\,\phi_4,\hspace{3mm} S_2=400\,\phi_6-21\,\phi_4,\hspace{3mm} S_3=40000\,\phi_6-2877\,\phi_4.\\
\end{array}
\label{atsh1}
\ee
The region of absolute stability is drawn in Figure~\ref{stab1}. The expressions for the phase-lag and dissipation associated to this method are given by
\[
\phi(\nu,z)=\ds\frac{23\,\epsilon}{378000\,(\omega^2+\epsilon)}\,H^7+{\mathcal O}(H^9),\qquad d(\nu,z)=-\ds\frac{37\,\epsilon}{216000(\omega^2+\epsilon)}\,H^6+{\mathcal O}(H^8).
\]

\begin{figure}[htp!]
\begin{center}
\begin{tabular}{c}
\includegraphics[width=12 cm]{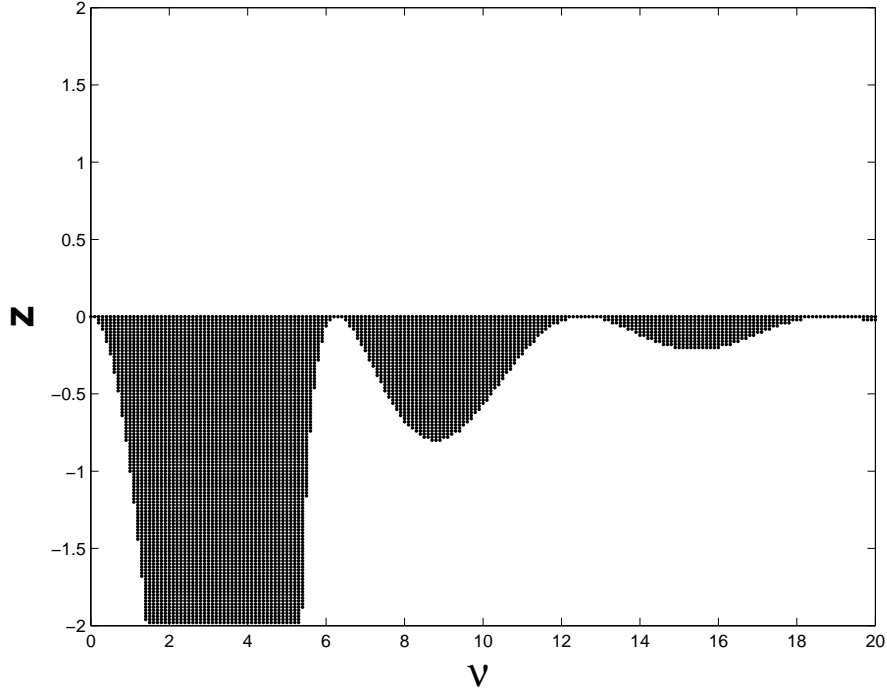}\\
\end{tabular}
\caption{\label{stab2}$\nu-z$ plot for ATSH method~(\ref{atsh2}).}
\end{center}
\end{figure}

{\bf * ATSH Method with phase-lag order eight}
\newline
Following Corollary~\ref{corr1} the condition that imposes phase-lag order eight is the same as that for the classical method. In the classical case, phase-lag order eight is achieved when $c_3=25/28$, see Franco~(2006a). Guided by Franco's method, we conclude with the coefficients
\be
\begin{array}{l}
a_{31}=\ds\frac{1325}{43904},\qquad a_{32}=\ds\frac{35775}{43904},\qquad a_{41}=\ds\frac{28\,S_1\,S_2\,(18816\,\phi_6^2-2186\,\phi_6\,\phi_4+53\,\phi_4^2)}{4293\,\phi_4^4},\\[5mm]

a_{42}=-\ds\frac{S_1\,S_2\,(526848\,\phi_6^2-51800\,\phi_6\,\phi_4+475\,\phi_4^2)}{2025\,\phi_4^4},\qquad a_{43}=\ds\frac{1568\,S_1\,S_2\,S_3\,\phi_6}{107325\,\phi_4^4},\\[5mm]

b_1=\ds\frac{2\,(9408\,\phi_6-625\,\phi_4)\,\phi_4}{53\,S_2},\\[5mm] b_2=\ds\frac{2\,(1418\,\phi_4^2-625\,\phi_4\,\phi_2-18816\,\phi_6\,\phi_4+8400\,\phi_2\,\phi_6)}{25\,S_1},\\[5mm]

b_3=\ds\frac{2458624\,(12\,\phi_6-\phi_4)\,\phi_4}{1325\,S_3},\qquad b_4=\ds\frac{162\,\phi_4^4}{S_1\,S_2\,S_3},\qquad c_3=\ds\frac{25}{28},\qquad c_4=\ds\frac{S_1}{3\,\phi_4},\\[5mm]
S_1=336\,\phi_6-25\,\phi_4,\qquad S_2=168\,\phi_6-11\,\phi_4,\qquad S_3=9408\,\phi_6-775\,\phi_4.
\end{array}
\label{atsh2}
\ee
The region of absolute stability is drawn in Figure~\ref{stab2}. The phase-lag and dissipation for this method are
\[
\phi(\nu,z)=-\ds\frac{(199\,\omega^2+182\,\epsilon)\,\epsilon}{101606400\,(\omega^2+\epsilon)^2}\,H^9+{\mathcal O}(H^{11}),\qquad d(\nu,z)=-\ds\frac{\epsilon}{20160(\omega^2+\epsilon)}\,H^6+{\mathcal O}(H^8).
\]

\subsubsection{Zero-dissipative fourth-order method with phase-lag order six}
\begin{figure}[htp!]
\begin{center}
\begin{tabular}{c}
\includegraphics[width=12 cm]{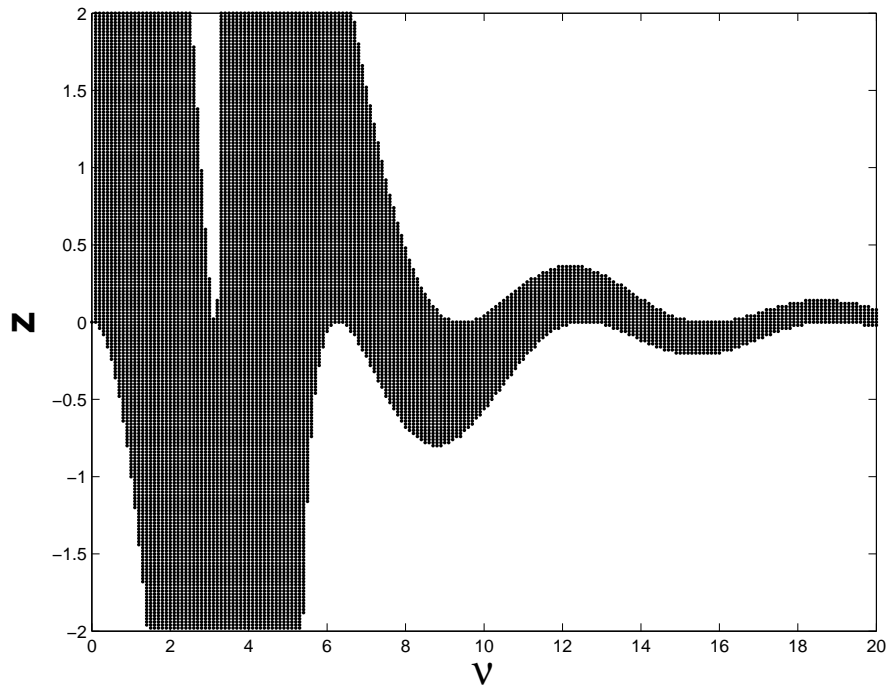}\\
\end{tabular}
\caption{\label{stab3}$\nu-z$ plot for ATSH method~(\ref{atsh3}).}
\end{center}
\end{figure}

Here we investigate how we can obtain zero-dissipative methods.  Following Theorem~\ref{phaseorderatsh} the method has phase-lag order six when
\be
b^T\,A^2\,c=0,\qquad b^T\,A^2\,e=2\,\phi_6(\nu).
\label{phaselag6}
\ee
We find $c_3=1$ which is incompatible with the fifth-order conditions~(\ref{ordercond5}), and the algebraic order of the method should be restricted to four. Solving equations (\ref{simplcond}), (\ref{ordercond4}) and (\ref{phaselag6}) we obtain the coefficients in terms of arbitrary parameters $c_3$ and $c_4$. The error constant $E^{ATSH}_5$~(\ref{errca}) should be as small as possible so we have that $c_4=(5\,c_3-2)/(5\,c_3-5)$, just like Franco's original case. It is easy to verify that the method reaches order five for linear systems of ODEs
\be
y''=-\omega^2\,y+g(x).
\label{linprob}
\ee
In the classical case the free parameter $c_3$ is chosen so that the resulting method is optimized for the class of linear problems~(\ref{linprob}). Here, in order to calculate the error constant when solving~(\ref{linprob}), we have to consider the coefficients of the $7$th-order elementary differentials $f^{(5)}(x)(y',y',y',y',y')$, $f^{(1)}(y)\left(f^{(3)}(x)(y',y',y')\right)$ and $\omega^2\,f^{(3)}(x)(y',y',y')$. The other $7$th-order elementary differentials remain zero for (\ref{linprob}). Minimizing this error constant we obtain $c_3=13/20$. For comparison, in the classical case Franco~(2006a) obtained $c_3=33/50$. The following coefficients are found 
\be
\begin{array}{l}
a_{31}=0,\qquad a_{32}=\ds\frac{429}{800},\qquad a_{41} = \ds\frac{38200\,\phi_6}{79233\,\phi_4},\qquad  a_{42}=-\ds\frac{5\,(7640\,\phi_6+637\,\phi_4)}{31213\,\phi_4},\\[5mm]
a_{43}=\ds\frac{764000\,\phi_6}{1030029\,\phi_4},\qquad b_1=-\ds\frac{6\,\phi_4}{11},\qquad b_2=-\ds\frac{596\,\phi_4}{65}+2\,\phi_2,\qquad b_3=\ds\frac{128000\,\phi_4}{27313},\\[5mm]
b_4=\ds\frac{4802\,\phi_4}{955},\qquad c_3=\ds\frac{13}{20},\qquad c_4=-\ds\frac{5}{7}.\\[5mm]
\end{array}
\label{atsh3}
\ee
The region of periodicity is drawn in Figure~\ref{stab3}, and the phase-lag is
\[
\phi(\nu,z)=-\ds\frac{\epsilon}{40320\,(\omega^2+\epsilon)}\,H^7+{\mathcal O}(H^{9}).
\]

\section{Numerical experiments}
\label{secnumexp}
\begin{figure}[ht!]
\begin{center}
\begin{tabular}{cc}
\includegraphics[width=12 cm]{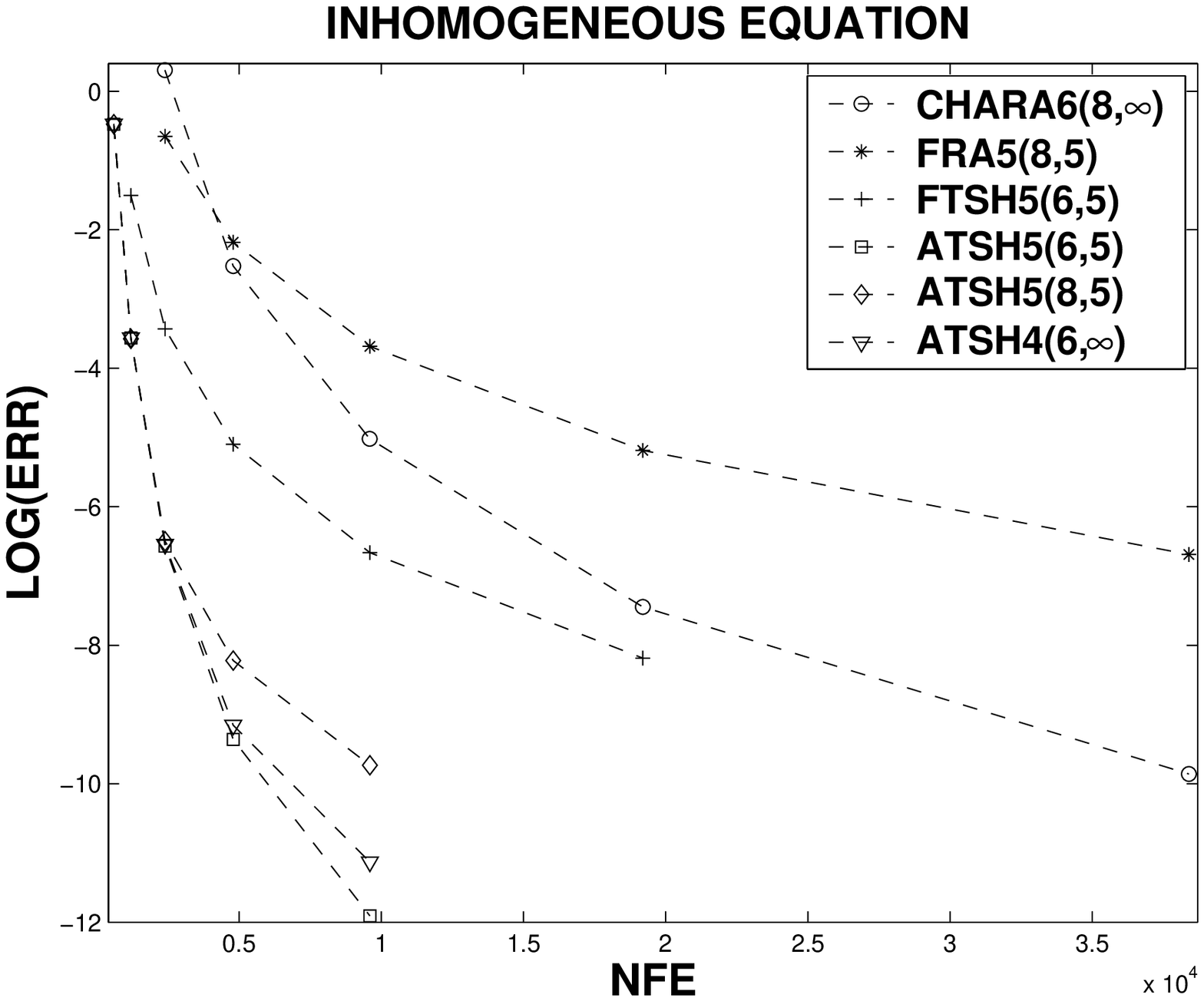}\\
\includegraphics[width=12 cm]{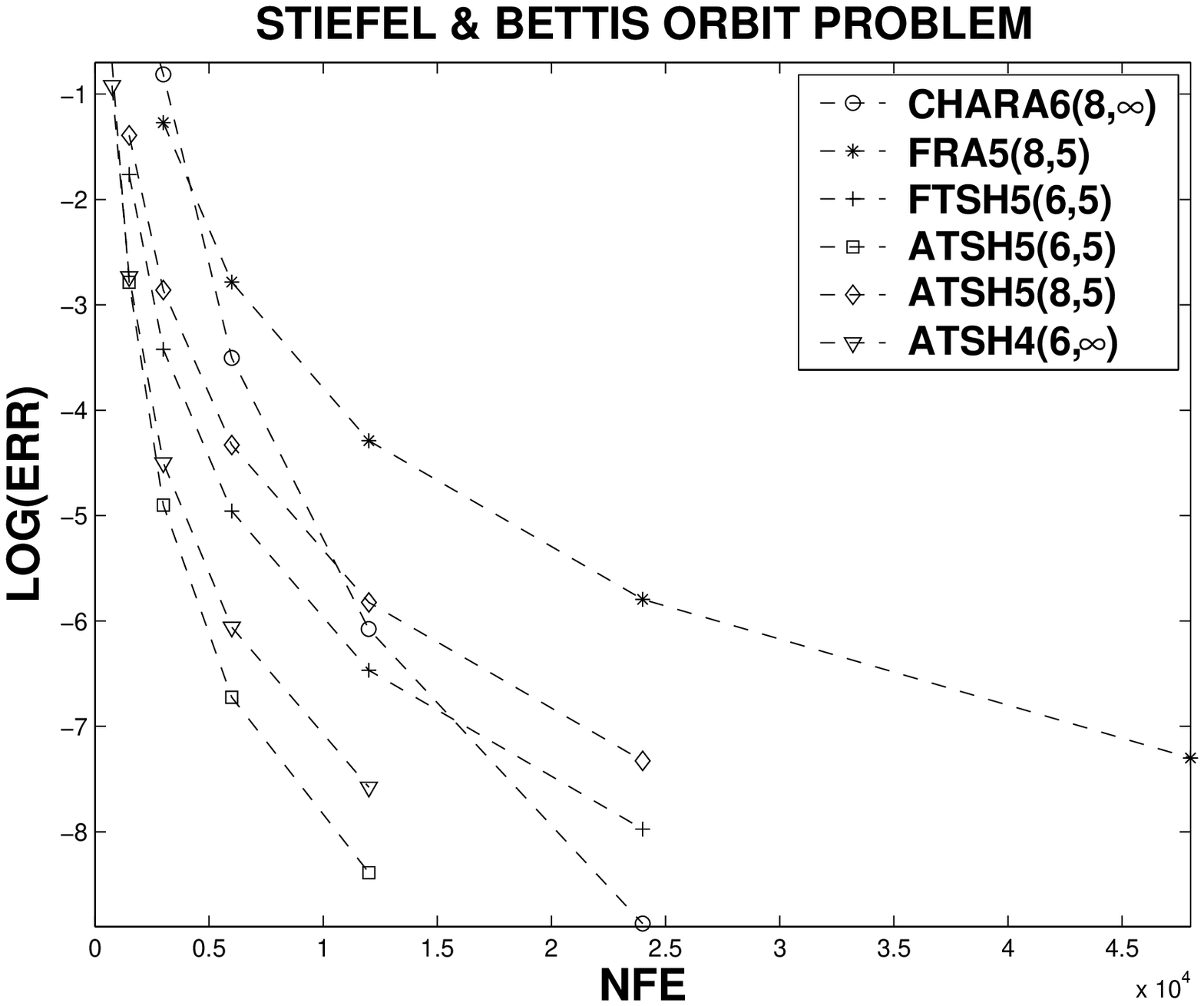}\\
\end{tabular}
\caption{\label{eff1}Efficiency curves of the methods for Problems 1--2.}
\end{center}
\end{figure}

\begin{figure}[ht!]
\begin{center}
\begin{tabular}{cc}
\includegraphics[width=12 cm]{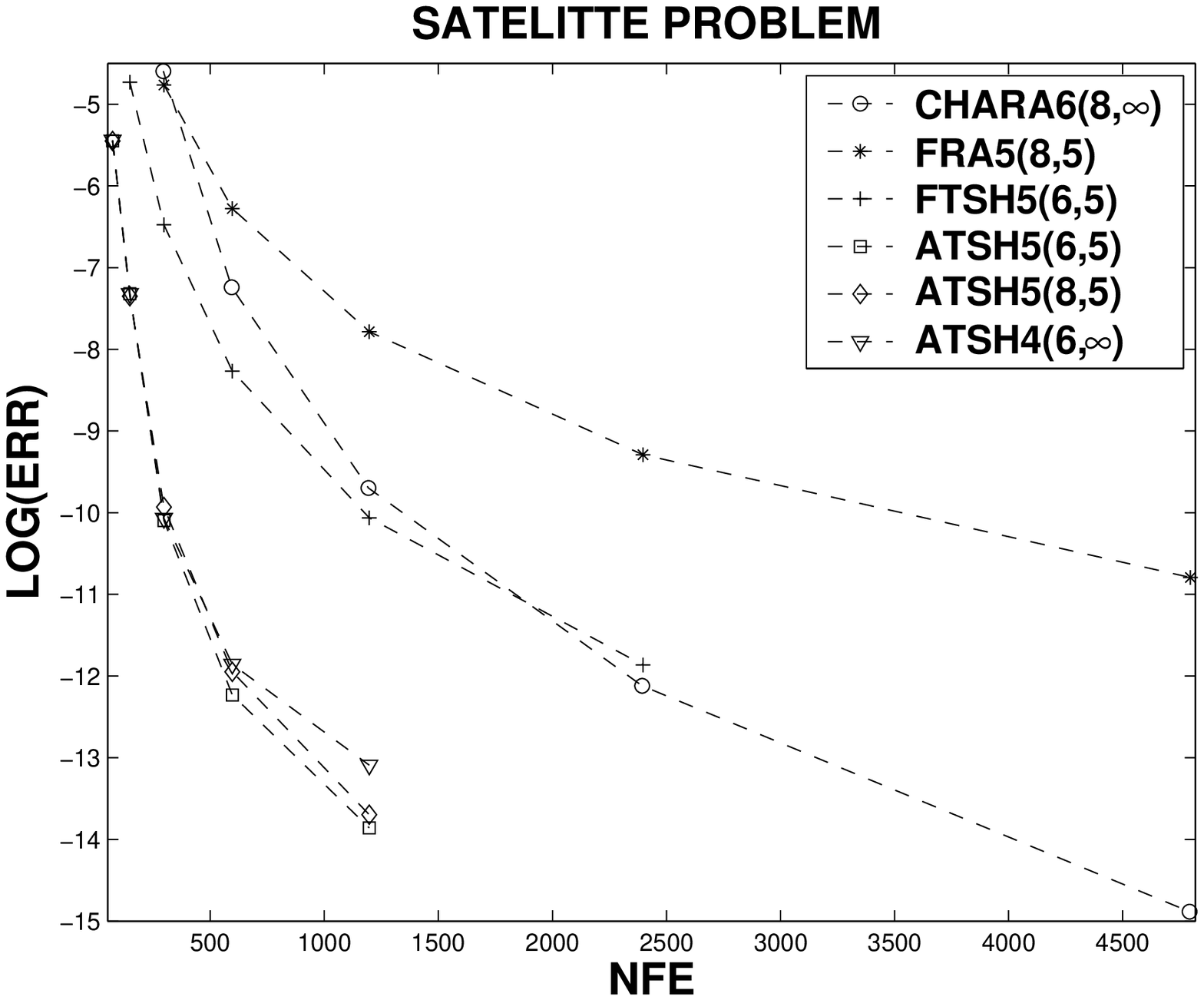}\\
\includegraphics[width=12 cm]{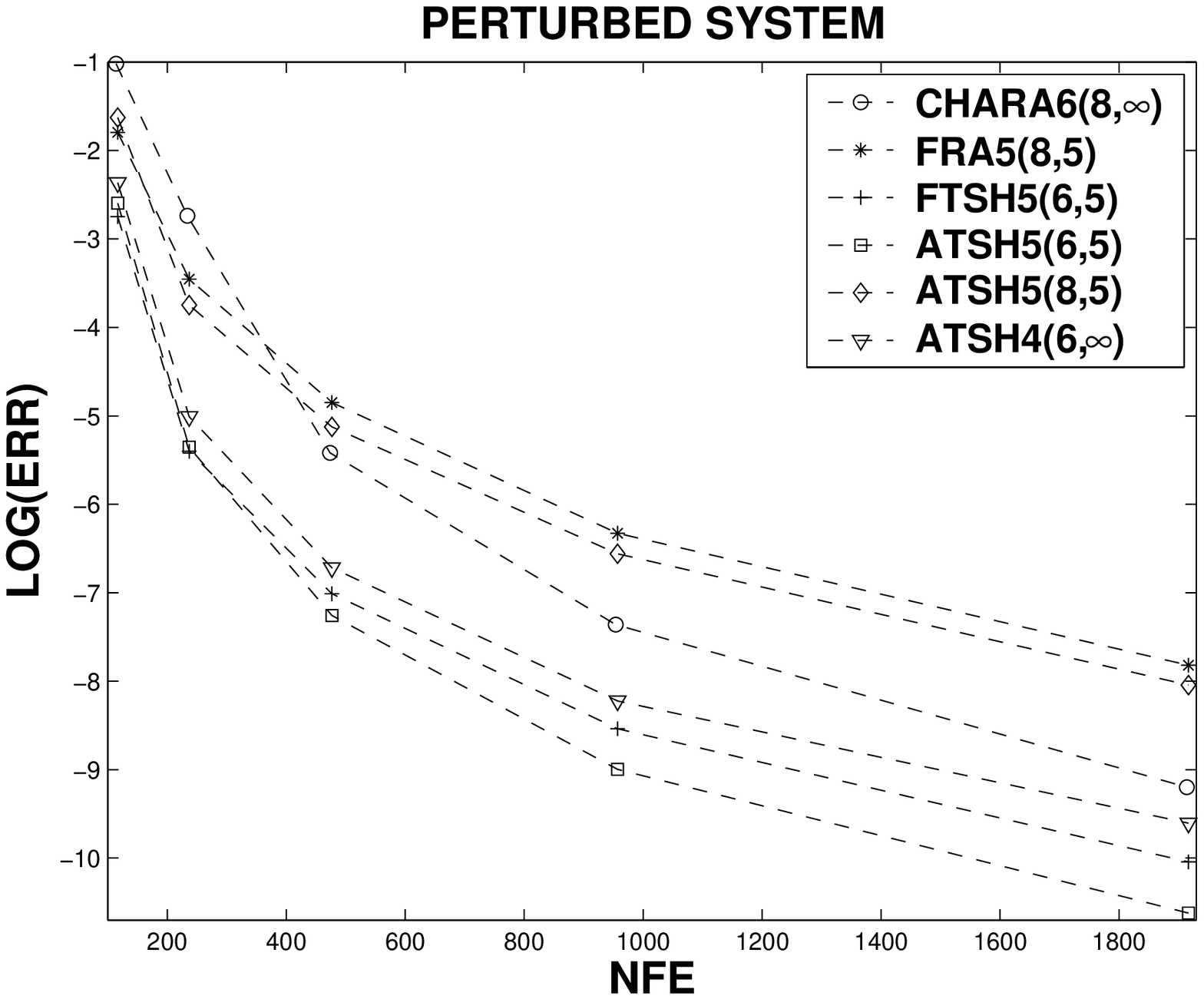}\\ 
\end{tabular}
\caption{\label{eff2}Efficiency curves of the methods for Problems 3--4.}
\end{center}
\end{figure}

In order to evaluate the effectiveness of the new method derived above we consider several model problems. The new method have been compared with other explicit TSH codes proposed in the literature. The criterion used in the numerical comparisons is the usual test based on computing the maximum global error over the whole integration interval. In Figures~\ref{eff1}--\ref{eff2} we have depicted the efficiency curves for the tested codes. These figures show the decimal logarithm of the maximum global error versus the computational error measured by the number of function evaluations required by each code. The algorithms used in the comparisons have been denoted by
\begin{itemize}
\item CHARA6(8,$\infty$): Zero-dissipative method derived by Chawla \& Rao~(1987).
\item FRA5(8,5): Classical method derived by Franco~(2006a).
\item FTSH5(6,5): Phase-fitted and amplification-fitted method derived by Van de Vyver~(2006).
\item ATSH5(6,5): ATSH method~(\ref{atsh1}).
\item ATSH5(8,5): ATSH method~(\ref{atsh2}).
\item ATSH4(6,$\infty$): ATSH method~(\ref{atsh3}).
\end{itemize}
Here, A(B,C) means that the method has algebraic order A, phase-lag order B and dissipation order C.

We have used the following five model problems:

{\bf Problem 1.} {\em An inhomogeneous equation studied by van der Houwen and~Sommeijer (1987)}
\[
y''=-100\,y+99\,\sin(x),\qquad y(0)=1,\qquad y'(0)=11.
\]
The exact solution is given by:
\[
y(x)=\cos(10\,x)+\sin(10\,x)+\sin(x). 
\]
It consists of a rapidly and slowly oscillating function; the slowly varying function is due to the inhomogeneous term. The equation has been solved in the interval $[0,100]$ with fitted frequency is $\omega=10$. The numerical results stated in Fig.~\ref{eff1} have been computed with stepsizes $h=2^{-j}$, $j=2,\ldots,6$ for CHARA6(8,$\infty$) and FTSH5(6,5), $j=3,\ldots,7$ for FRA5(8,5) and $j=1,\ldots,5$ for ATSH5(6,5), ATSH5(8,5) and ATSH4(6,$\infty$).

{\bf Problem 2.} {\em An ``almost periodic'' orbit problem studied by Stiefel and Bettis~(1969)}
\[
z''=-z+0.001\,e^{i\,x},\qquad z(0)=1,\qquad z'(0)=0.9995\,i.
\]
The equation has been solved in the interval $[0,1000]$ with fitted frequency $\omega=1$.
The exact solution is given by:
\[
z(x)=(1-0.0005\,i\,x)\,e^{i\,x}.
\]
The solution represents a motion of a perturbation of a circular orbit in the complex plane. The problem may be solved either as a single equation in complex arithmetic or as a pair of uncoupled equations. The numerical results stated in Fig.~\ref{eff1} have been computed with stepsizes $h=2^{-j}$, $j=-2,\ldots,2$ for CHARA6(8,$\infty$), ATSH5(6,5) and ATSH4(6,$\infty$), $j=-1,\ldots,3$ for FTSH5(6,5) and ATSH5(8,5), $j=0,\ldots,4$ for FRA5(8,5).

{\bf Problem 3.} {\em A satellite problem studied by Ferr\'andiz et al.~(1992)}
\newline
We consider the problem of determining the position of an earth satellite. The equations of motion have been expressed in focal variables (see Ferr\'andiz (1988) and Ferr\'andiz~{\em et al.}~(1992)). The coordinates of the basic set of focal variables are three components $(y_1,y_2,y_3)$ of the direction vector of the particle and the inverse $u$ of the radial distance. In this formulation the satellite problem can be formulated in four decoupled pertubed harmonic oscillators with unit frequency:
\be
\begin{array}{lll}
y_i''+y_i&=&Q_i,\qquad i=1,2,3,\\\\
u''+u&=&\ds\frac{\mu}{c^2}+Q,
\end{array}
\label{sat1}
\ee
where $\mu$ is the reduced mass, while $Q_i$ and $Q$ denote the corresponding perturbation terms. We consider the almost periodic equatorial orbit with the zonal harmonic coefficient $J_2$ taken as the perturbation parameter. We have neglected higher order terms of $J_2$. The system of equations (\ref{sat1}) can be written in the form
\be
\begin{array}{lll}
y_i''+y_i&=&0,\qquad i=1,2,3,\\\\
u''+u&=&\ds\frac{\mu}{c^2}+12\,\ds\frac{J_2}{c^2}u^2,
\end{array}
\label{sat2}
\ee
where $c$ is the angular momentum and it can be considered as a constant. The solutions of the first three oscillators are trivial thus we are focused on the last equation. We consider the domain of integration $[\pi,100]$. The initial conditions are given by
\[
u(\pi)=\ds\frac{\mu\,(1-e)}{c^2},\qquad u'(\pi)=0.
\]
For our numerical purpose we consider orbits with eccentricity $e=0.99$. In this case: 
\[ 
\begin{array}{ll}
\ds\frac{\mu}{c^2}=\ds\frac{100}{20895},& \ds\frac{J_2}{c^2}=\ds\frac{50}{20895000}.\\
\end{array}
\]
The error has been calculated using a reference solution obtained by means of the perturbation techniques developed by Farto~{\em et al.}~(1998). The numerical results stated in Fig.~\ref{eff2} have been computed with stepsizes $h=(1-\pi/100)\,2^{-j}$, $j=-1,\ldots,3$ for CHARA6(8,$\infty$) and FTSH5(6,5), $j=0,\ldots,4$ for FRA5(8,5), $j=-2,\ldots,2$ for ATSH5(6,5), ATSH5(8,5) and ATSH4(6,$\infty$).

{\bf Problem 4.} {\em A perturbed system studied by Franco (2002)}
\newline
As an example of a system we consider
\[
\begin{array}{l}
y_1''=-25\,y_1-\epsilon\,(y_1^2+y_2^2)+\epsilon\,f_1(x),\qquad y_1(0)=1,\qquad y'_1(0)=0,\\
y_2''=-25\,y_2-\epsilon\,(y_1^2+y_2^2)+\epsilon\,f_2(x),\qquad y_2(0)=\epsilon,\qquad y'_2(0)=5,\\
\end{array}
\]
where
\[
\begin{array}{lll}
f_1(x)&=&1+\epsilon^2+2\,\epsilon\,\sin(5\,x+x^2)+2\,\cos(x^2)+(25-4\,x^2)\,\sin(x^2),\\
f_2(x)&=&1+\epsilon^2+2\,\epsilon\,\sin(5\,x+x^2)-2\,\sin(x^2)+(25-4\,x^2)\,\cos(x^2).\\
\end{array}
\]
In our test we choose $\epsilon=10^{-3}$. The system has been solved in the interval $[0,5]$ with $\omega=5$. The analytical solution is given by:
\[
y_1(x)=\cos(5\,x)+\epsilon\,\sin(x^2),\qquad y_2(x)=\sin(5\,x)+\epsilon\,\cos(x^2).
\]
The numerical results stated in Fig.~\ref{eff2} have been computed with stepsizes $h=2^{-j}$, $j=1,\ldots,5$ for CHARA6(8,$\infty$), $j=2,\ldots,6$ for the other codes.

\section{Conclusions}
Scheifele's $G$-functions methods are designed in such a way that the exact integration of the homogeneous solution of perturbed oscillators~(\ref{pert}) is automatically included. The methods take care with the evaluation of the inhomogeneous part of (\ref{pert}), i.e. $g(x,y)$. We have applied Scheifele's approach to TSH methods for an accurate and efficient integration of~(\ref{pert}). The resulting methods, called ATSH methods, have coefficients dependent on $\nu=\omega\,h$, where $\omega$ is a specified angular frequency. Classical TSH methods are the limiting forms of ATSH methods as $\nu\rightarrow 0$. 

This paper provides a theoretical framework for the derivation of ATSH methods. One of our main aims is to develop the order conditions for this new type of methods. It is found that ATSH methods share some important properties with the corresponding classical TSH methods such as zero-stability, the dissipation order and, under some conditions, with the phase-lag order. On the contrary, the stability properties are very different from the classical method and they depend on the fitted frequency and the stepsize. When the main frequency of the problem is exactly known stability problems will never occur, except for a discrete set of exceptional values of the stepsize. When the dominant frequency is not exactly known some care is required when selecting the stepsize.

In particular, we have demonstrated the validity of the theory with explicit fourth- and fifth-order ATSH methods. The new methods are adaptations of the classical TSH methods of Franco~(2006a). In most cases, the dissipative ATSH method (\ref{atsh1}) with minimized error constant outperforms all the other methods considered. In contrast with the results of the phase-fitted and amplification-fitted methods of Van de Vyver (2007b), it turns out that the accuracy of ATSH methods is mostly determined by its usual local truncation error rather than by its phase-lag.

Our task is restricted to scalar equations or systems involving only one frequency. When solving systems with more than one frequency, or more general, systems of the form
\be
y''=K\,y+g(x,y),
\label{sys}
\ee
the resulting methods have coefficients which are functions of the matrix $h^2\,K$. So their evaluation is not direct. To overcome this difficulty, together with some other troubles, Franco (2006b) has modified ARKN methods for oscillatory systems of the form~(\ref{sys}). The extension of Franco's approach to the ATSH methods considered here might be an interesting suggestion for some future work.
\section*{Acknowledgments}
This research was supported by ``Grant 0T/04/21 of Onderzoeksfonds K.U. Leuven" and ``Scholarship BDB-B/05/06 of K.U. Leuven".
\begin{center}
{\large REFERENCES}
\end{center}
CHAWLA, M.~M.~(1984) Numerov made explicit has better stability. {\em BIT}, {\bf 24}, 117-118. 
CHAWLA, M.~M.~\& RAO, P.~S.~(1987) An explicit sixth-order method with phase-lag of order eight for $y'' = f(t,y)$. {\em J. Comput. Appl. Math.}, {\bf 17}, 365--368.
\newline
COLEMAN, J.~P.~(1989) Numerical methods for $y''=f(x,y)$ via rational approximations  for the cosine. {\em IMA J. Numer. Anal.}, {\bf 9}, 145--165.
\newline
COLEMAN, J.~P.~(2003) Order conditions for a class of two-step methods for $y''=f(x,y)$. {\em IMA J. Numer. Anal.}, {\bf 23}, 197--220.
\newline
COLEMAN, J.~P. \& IXARU, L.~GR (1996) P-stability and exponential-fitting methods for $y''=f(x, y)$. {\em IMA J. Numer. Anal.}, {\bf 16}, 179-199.
\newline
FAIR\'EN, V., MART\'IN, P. \& FERR\'ANDIZ, J.~M. (1994) Numerical tracking of small deviations from analytically known periodic orbits. {\em Computers in Physics}, {\bf 8}, 455--461.
\newline
FARTO, J.~M., GONZ\'ALEZ, A.~B. \& MART\'IN, P. (1998) An algorithm for the systematic construction of solutions to perturbed problems. {\em Comput. Phys. Commun.},~{\bf 111}, 110--132.
\newline
FERR\'ANDIZ, J.~M. (1988) A general canonical transformation increasing the number of variables with applications to the two-body problem. {\em Celest. Mech.}, {\bf 41}, 343--357.
\newline
FERR\'ANDIZ, J.~M., SANSATURIO, M.~E.~\& POJMAN, J.~R.~(1992) Increased accuracy of computations in the main satellite problem through linearization methods. {\em Celest. Mech. Dynam. Astronom.}, {\bf 53}, 347--363.
\newline
FRANCO, J.~M.~(2002) Runge-Kutta-Nystr\"om methods adapted to the numerical integration of perturbed oscillators. {\em Comput. Phys. Commun.}, {\bf 147}, 770--787.
\newline
FRANCO, J.~M.~(2005) Stability of explicit ARKN methods for perturbed oscillators. {\em J. Comput. Appl. Math.}, {\bf 173}, 389--396. 
\newline
FRANCO, J.~M.~(2006a) A class of explicit two-step hybrid methods for second-order IVPs. {\em J. Comput. Appl. Math.}, {\bf 187}, 41--57.
\newline
FRANCO, J.~M.~(2006b) New methods for oscillatory systems based on ARKN methods. {\em Appl. Numer. Math.}, {\bf 56}, 1040-1053.  \newline
GONZ\'ALEZ, A.~B., MART\'IN, P. \& FARTO, J.~M. (1999) A new family of Runge-Kutta type methods for the numerical integration of perturbed oscillators. {\em Numer. Math.}, {\bf 82}, 635--646.
\newline
HAIRER, E., N\O RSETT, S.~P. \& WANNER, G. (1993) {\em Solving Ordinary Differential Equations~I, Nonstiff Problems}, 2nd edn. Springer Series in Computational Mathematics, Berlin.
\newline
HENRICI, P. (1962) {\em Discrete Variable Methods in Ordinary Differential Equations}. Wiley, New York.
\newline
IXARU, L.~GR. \& RIZEA, M. (1987) Numerov method maximally adapted to the Schr\"odinger equation. {\em J. Comput. Phys.}, {\bf 73}, 306--324.
\newline
LAMBERT, J.~D. \& WATSON, I.~A. (1976) Symmetric multistep methods for periodic initial-value problems. {\em J. Inst. Math. Appl.}, {\bf 18}, 189--202.
\newline
L\'OPEZ, D.~J., MART\'IN, P. \& FARTO, J.~M. (1999) Generalization of St\"ormer method for perturbed oscillators without explicit first dervatives. {\em J. Comput. Appl. Math.}, {\bf 111}, 123--132.
\newline
MART\'IN, P. \& FERR\'ANDIZ, J.~M.~(1997) Multistep numerical methods based on the Scheifele $G$-functions with application to satellite dynamics. {\em SIAM J. Numer. Anal.}, {\bf 34}, 359--375.
\newline
PETZOLD, L.~R., JAY, L.~O \& YEN, J. (1997) Numerical solution of highly oscillatory ordinary differential equations. {\em Numerica Acta}, 437--483.
\newline
SCHEIFELE, G.~(1971) On the numerical integration of perturbed linear oscillating systems. {\em Z. Angew. Math. Phys.}, {\bf 22}, 186--210. 
\newline
SIMOS, T.~E. (1999) Explicit eight order methods for the numerical integration of initial-value problems with periodic or oscillating solutions. {\em Comput. Phys. Commun.}, {\bf 119}, 32--44. 
\newline
STIEFEL, E. \& BETTIS, D.~G.~(1969) Stabilization of Cowell's method. {\em Numer. Math.}, {\bf 13}, 154--175.
\newline
TSITOURAS, CH.~(2003) Explicit Numerov type methods with reduced number of stages. {\em Comput. Math. Appl.}, {\bf 45}, 37--42.
\newline
VAN DE VYVER, H. (2006) A phase-fitted and amplification-fitted two-step hybrid method for second-order periodic initial value problems. {\em Internat. J. Modern Phys.~C}, {\bf 17}, 663--675.
\newline
VAN DE VYVER, H. (2007a) An adapted explicit hybrid method of Numerov type for the numerical integration of perturbed oscillators. {\em Appl. Math. Comp.}, corrected proof available online via ScienceDirect.
\newline
VAN DE VYVER, H. (2007b) Phase-fitted and amplification-fitted two-step hybrid methods for $y''=f(x,y)$. {\em J. Comput. Appl.  Math.}, corrected proof available online via ScienceDirect.
\newline
VAN DER HOUWEN, P.~J. \& SOMMEIJER, B.~P.~(1987) Explicit Runge-Kutta\newline(-Nystr\"om) methods with reduced phase errors for computing oscillating solutions. {\em SIAM J. Numer. Anal.}, {\bf 24}, 595--617.
\newline
VIGO-AGUIAR, J., SIMOS, T.~E. \& FERR\'ANDIZ, J.~M. (2004) Controlling the error growth in long-term numerical integration of perturbed oscillations in one or several frequencies. {\em Proc. R. Soc. Lond.} A, {\bf 460}, 561--567.

\end{document}